\documentclass{amsart}
\DeclareMathSymbol{\twoheadrightarrow}
{\mathrel}{AMSa}{"10}

\def\Q{{\mathbf Q}}
\def\Z{{\mathbf Z}}

\def\R{{\mathfrak  R}}
\def\F{{\mathbf F}}

\def\U{{\mathbf U}}
\def\PGU{\mathrm{PGU}}
\def\Sz{\mathbf{Sz}}

\def\Sn{{\mathbf S}_n}
\def\An{{\mathbf A}_n}

\def\Gal{\mathrm{Gal}}
\def\Perm{\mathrm{Perm}}
\def\tr{\mathrm{tr}}
\def\diag{\mathrm{diag}}

\def\End{\mathrm{End}}
\def\Aut{\mathrm{Aut}}

\def\I{\mathrm{Id}}

\def\fchar{\mathrm{char}}

\def\GL{\mathrm{GL}}
\def\SU{\mathrm{SU}}
\def\St{\mathrm{St}}
\def\PSU{\mathrm{PSU}}
\def\SL{\mathrm{SL}}

\def\PSL{\mathrm{PSL}}

\def\dim{\mathrm{dim}}

\def\P{{\mathbf P}}

\def\f{{\mathcal F}}

\newtheorem{thm}{Theorem}[section]
\newtheorem{lem}[thm]{Lemma}

\theoremstyle{definition}
\newtheorem{defn}[thm]{Definition}

\newtheorem{rem}[thm]{Remark}
\newtheorem{rems}[thm]{Remarks}
\title{Hyperelliptic jacobians and $\U_3(2^m)$}
\author[Yuri G. Zarhin]{Yuri G. Zarhin}
\address{Department of Mathematics, Pennsylvania State University,
University Park, PA 16802, USA}
\email{zarhin\char`\@math.psu.edu}
\thanks{Partially supported by NSF grant DMS-0070664}

\begin{document}
\maketitle
\section{Introduction}

In \cite{Zarhin} the author  proved that
in characteristic $0$ the
jacobian $J(C)=J(C_f)$  of a hyperelliptic curve
$$C=C_f:y^2=f(x)$$
has only trivial endomorphisms
over an algebraic closure $K_a$ of the ground field $K$
if the Galois group $\Gal(f)$ of the irreducible polynomial
$f \in K[x]$ is ``very big". Namely, if $n=\deg(f) \ge 5$
and $\Gal(f)$ is either the symmetric group $\Sn$ or the alternating group $\An$
 then the ring $\End(J(C_f))$ of $K_a$-endomorphisms of $J(C_f)$ coincides with $\Z$.  Later the author \cite{Zarhin2} proved that $\End(J(C_f))=\Z$ for an infinite series
of $\Gal(f)=\PSL_2(\F_{2^r})$ and $n=2^{r}+1$
(with $\dim(J(C_f))=2^{r-1}$) or when $\Gal(f)$ is
 the Suzuki group $\Sz(2^{2r+1})$ and $n=2^{2(2r+1)}+1$
(with $\dim(J(C_f))=2^{4r+1}$). We refer the reader to
\cite{Mori1}, \cite{Mori2}, \cite{Katz1}, \cite{Katz2},
\cite{Masser}, \cite{Zarhin}, \cite{Zarhin2}, \cite{ZarhinMRL2}
for a discussion of known results about, and examples of,
hyperelliptic jacobians without complex multiplication.

We write $\R=\R_f$ for the set of roots of $f$ and consider $\Gal(f)$ as the corresponding permutation group of $\R$.
Suppose $q=2^m>2$ is an integral power of $2$ and $\F_{q^2}$ is a finite field consisting of $q^2$ elements. Let us consider a non-degenerate Hermitian (wrt $x\mapsto x^q$) sesquilinear form on $\F_{q^2}^3$.
In the present paper we prove that
$\End(J(C_f))=\Z$ when  $\R_f$ can be identified with the corresponding ``Hermitian curve" of isotropic lines in the projective plane
$\P^2(\F_{q^2})$ in such a way  that $\Gal(f)$,  becomes either the projective unitary group $\PGU_3(\F_{q})$ or the projective special unitary group $\U_3(q):=\PSU_3(\F_q)$. In this case $n=\deg(f)=q^3+1=2^{3m}+1$ and $\dim(J(C_f))=q^3/2=2^{3m-1}$.

Our proof is based on an observation that the Steinberg
representation is the only absolutely irreducible nontrivial
representation (up to an isomorphism) over $\F_2$ of $\U_3(2^m)$,
whose dimension is a power of $2$.

\section{Main results}
\label{mainr}
Throughout this paper we assume that $K$ is a field with $\fchar(K) \ne 2$. We fix its algebraic closure $K_a$ and
 write $\Gal(K)$ for the absolute Galois group $\Aut(K_a/K)$. If $X$ is an abelian variety defined over $K$
 then we write $\End(X)$ for the ring of  $K_a$-endomorphisms of $X$.

Suppose $f(x) \in K[x]$ is a separable polynomial
 of degree $n\ge 5$. Let $\R=\R_f \subset K_a$ be
 the set of roots of $f$, let $K(\R_f)=K(\R)$ be
 the splitting field of $f$ and $\Gal(f):=\Gal(K(\R)/K)$ the Galois group of $f$, viewed as a subgroup of $\Perm(\R)$.
Let $C_f$ be the hyperelliptic curve $y^2=f(x)$. Let  $J(C_f)$ be
its jacobian, $\End(J(C_f))$ the ring of $K_a$-endomorphisms of $J(C_f)$.

\begin{thm}
\label{main2}
Assume that there exist a positive integer $m>1$ such that $n=2^{3m+1}+1$ and $\Gal(f)$ contains
 a subgroup isomorphic  to $\U_3(2^m)$.
Then either $\End(J(C_f))=\Z$ or  $\fchar(K)>0$ and
$J(C_f)$ is a supersingular abelian variety.
\end{thm}

We will prove Theorems \ref{main2}  in \S \ref{final}.

\section{Permutation groups, permutation modules and very simplicity}
\label{permute}

Let $B$ be a finite set consisting of $n \ge 5$ elements. We write $\Perm(B)$ for the group of permutations of $B$. A choice of ordering on $B$ gives rise to an isomorphism
$$\Perm(B) \cong \Sn.$$
Let $G$ be a  subgroup of $\Perm(B)$.
For each $b \in B$ we write $G_b$ for the stabilizer of $b$ in $G$; it is a subgroup of $G$.
Further we always assume that $n$ is {\sl odd}.

\begin{rem}
\label{transitive}
Assume that the action of $G$ on $B$ is transitive.
It is well-known that each $G_b$ is  of index $n$ in $G$ and
 all the $G_b$'s are conjugate in $G$.
 Each conjugate of $G_b$ in $G$ is the stabilizer of a point of $B$.
 In addition, one may identify the $G$-set $B$ with the set of cosets $G/G_b$
with the standard action by $G$.
\end{rem}

 We write $\F_2^B$
for the $n$-dimensional $\F$-vector space of maps $h:B \to \F_2$.
The space $\F_2^B$ is  provided with a natural action of $\Perm(B)$ defined
as follows. Each $s \in \Perm(B)$ sends a map
 $h:B\to \F_2$ into  $sh:b \mapsto h(s^{-1}(b))$. The permutation module $\F_2^B$
contains the $\Perm(B)$-stable hyperplane
$$Q_B:=
\{h:B\to\F_2\mid\sum_{b\in B}h(b)=0\}$$
and the $\Perm(B)$-invariant line $\F \cdot 1_B$ where $1_B$ is the constant function $1$.
Since $n$ is odd,  there is a $\Perm(B)$-invariant splitting
$$\F_2^B=Q_B \oplus \F_2 \cdot 1_B.$$
Clearly,
$$\dim_{\F_2}(Q_B)=n-1$$
and $\F_2^B$ and $Q_B$  carry natural structures of $G$-modules.  Clearly, $Q_B$ is a faithful $G$-module.
It is also clear that the $G$-module $Q_B$ can be viewed as the reduction modulo $2$ of the
$\Q[G]$-module
$$(\Q^B)^0:=\{h:B\to\Q\mid\sum_{b\in B}h(b)=0\}.$$
It is well-known that the $\Q[G]$-module
$(\Q^B)^0$ is absolutely simple if and only if the action of $G$ on $B$ is doubly transitive
(\cite{SerreRep}, Sect. 2.3, Ex. 2).

\begin{rem}
\label{St2}
Assume that $G$ acts on $B$ doubly transitively
 and $\#(B)-1=\dim_{\Q}((\Q_B)^0)$ coincides with the largest power of $2$ dividing $\#(G)$.
 Then it follows from a theorem of Brauer-Nesbitt
 (\cite{SerreRep}, Sect. 16.4, pp. 136--137 ; \cite{Hump}, p. 249) that
 $Q_B$ is an absolutely simple $\F_2[G]$-module. In particular, $Q_B$ is  (the reduction
 of) the Steinberg representation \cite{Hump}, \cite{Curtis}.
\end{rem}

We refer to \cite{Zarhin2} for a discussion of the following definition.

\begin{defn}
Let $V$ be a vector space over a field $\F$, let $G$ be a group and
$\rho: G \to \Aut_{\F}(V)$ a linear representation of $G$ in $V$. We
say that the $G$-module $V$ is {\sl very simple} if it enjoys the
following property:

If $R \subset \End_{\F}(V)$ is an $\F$-subalgebra containing the
identity operator $\I$ such that

 $$\rho(\sigma) R \rho(\sigma)^{-1} \subset R \quad \forall \sigma \in G$$
 then either $R=\F\cdot \I$ or $R=\End_{\F}(V)$.
\end{defn}

\begin{rems}
\label{image}
\begin{enumerate}
\item[(i)]
If $G'$ is a subgroup of $G$ and the $G'$-module $V$ is very simple
then obviously the $G$-module $V$ is also very simple.

\item[(ii)]
A  very simple module is absolutely simple (see \cite{Zarhin2}, Remark 2.2(ii)).

\item[(iii)]
If $\dim_{\F}(V)=1$ then obviously the $G$-module $V$ is very simple.

\item[(iv)]
Assume that the $G$-module $V$ is  very simple and $\dim_{\F}(V)>1$.
Then  $V$ is not induced from a subgroup $G$ (except $G$ itself) and is not isomorphic to a tensor product
of two $G$-modules, whose $\F$-dimension is strictly less than $\dim_{\F}(V)$  (see \cite{Zarhin2}, Examples 7.1).

\item[(v)] If $\F=\F_2$ and $G$ is {\sl perfect} then  properties (ii)-(iv)
characterize the very simple $G$-modules (see \cite{Zarhin2}, Th. 7.7).
\end{enumerate}
\end{rems}

The following statement provides a criterion of very simplicity over $\F_2$.

\begin{thm}
\label{Very3}
Suppose  a positive integer $N>1$ and a group $H$ enjoy the following properties:
\begin{itemize}
\item
$H$ does not contain a subgroup of index dividing $N$ except $H$ itself.

\item
Let $N=ab$ be a factorization of $N$ into a product of two
positive integers $a>1$ and $b>1$. Then either
there does not exist an absolutely simple $\F_2[H]$-module of $\F_2$-dimension $a$ or
there does not exist an absolutely simple $\F_2[H]$-module of $\F_2$-dimension $b$.
\end{itemize}

Then each absolutely simple $\F_2[H]$-module of $\F_2$-dimension $N$ is very simple.
\end{thm}

\begin{proof}
This is Corollary 4.12 of \cite{Zarhin2}.
\end{proof}

\section{Steinberg representation}
\label{St}
We refer to \cite{Hump} and \cite{Curtis}  for a definition and basic properties of Steinberg representations.

Let us fix an algebraic closure of $\F_2$ and denote it by $\f$. We write $\phi: \f \to \f$ for the Frobenius automorphism
 $x \mapsto x^2$. Let $q=2^m$ be a positive integral power of two. Then
the subfield of invariants of $\phi^m:\f \to \f$ is a finite field $\F_q$ consisting of $q$ elements.
Let  $q'$ be an integral positive power of $q$.
 If $d$ is a positive integer and $i$ is a non-negative integer then for each
matrix $u \in \GL_d(\f)$ we write $u^{(i)}$ for the matrix obtained by raising each entry of $u$ to the $2^i$th power.

\begin{rem}
\label{prim}
Recall that an element $\alpha \in \F_q$ is called {\sl primitive} if $\alpha \ne 0$ and has multiplicative order $q-1$ in
 the cyclic multiplicative group $\F_q^*$.

Let $M<q-1$ be a positive integer. Clearly, the set
$$\mu_M(\F_q)=\{\alpha \in \F_q\mid \alpha^M=1\}$$
is a cyclic multiplicative subgroup of $\F_q^*$ and its order $M'$
divides both $M$ and $q-1$. Since $M<q-1$
and $q-1$ is odd, the ratio $(q-1)/M'$ is an {\sl odd} integer
$>1$. This implies that $3 \le (q-1)/M'$ and therefore
$$M'=\#(\mu_M(\F_q)) \le (q-1)/3.$$
\end{rem}

\begin{lem}
\label{trace}
Let $q>2$, let $d$ be a positive integer and let $G$ be a subgroup of $\GL_d(\F_{q'})$.
Assume  that one of the following two conditions holds:
\begin{enumerate}
\item[(i)]
There exists an element
$u \in G\subset\GL_d(\F_{q'})$, whose trace $\alpha$ lies in $\F_q^*$  and has multiplicative order $q-1$;
\item[(ii)]
There are exist a positive integer $r>\frac{q-1}{3}$, distinct
$\alpha_1, \cdots , \alpha_{r} \in \F_q^*$
 and elements
$$u_1, \cdots ,u_{r} \in G \subset \GL_d(\F_{q'})$$
such that the trace of $u_i$ is $\alpha_i$ for all $i=1,\cdots ,
r$.
\end{enumerate}
Let $V_0=\f^d$ and
$\rho_0: G \subset \GL_d(\F_{q'}) \subset \GL_d(\f)=\Aut_{\f}(V_0)$
be the natural $d$-dimensional representation of $G$ over $\f$. For each positive integer $i<m$
we define a $d$-dimensional $\f$-representation
$$\rho_i: G \to \Aut(V_i)$$
as the composition of
$$G \hookrightarrow \GL_d(\F_{q'}), \quad x \mapsto x^{(i)}$$
 and the inclusion map
$$\GL_d(\F_{q'}) \subset \GL_d(\f)\cong \Aut_{\f}(V_i).$$
Let $S$ be a subset of $\{0,1, \ldots m-1\}$.
Let us define a $d^{\#(S)}$-dimensional $\f$-representation $\rho_S$ of $G$ as the tensor product of
 representations $\rho_i$ for all $i \in S$. If $S$ is a proper subset of $\{0,1, \ldots m-1\}$
 then there exists an element $u \in G$ such that the trace of $\rho_S(u)$ does not belong to $\F_2$.
 In particular, $\rho_S$ could not be obtained by extension of scalars to $\f$ from a representation of $G$ over $\F_2$.
\end{lem}

\begin{proof}
Clearly,
$$\tr(\rho_i(u))=(\tr(\rho_0(u))^{2^i} \quad  \forall u \in G.$$
This implies easily that
$$\tr(\rho_S(u))=\prod_{i\in S}\tr(\rho_i(u))= (\tr(\rho_0(u))^M$$
 where $M=\sum_{i \in S}2^i$. Since $S$ is a {\sl proper} subset of $\{0,1, \ldots m-1\}$, we have
$$0<M < \sum_{i=0}^{m-1} 2^i=2^{m}-1=\#(\F_q^*).$$
Assume that the condition (i) holds.
Then there exists $u \in G$ such that $\alpha=\tr(\rho_0(u))$ lies in $\F_q^*$
and the exact multiplicative order of $\alpha$ is
 $q-1=2^m-1$.

This implies that $0 \ne \alpha^M \ne 1$. Since $\F_2=\{0,1\}$, we conclude that $\alpha^M \not\in \F_2$. Therefore
$$\tr(\rho_S(u))=(\tr(\rho_0(u))^M=\alpha^M \not\in\F_2.$$
Now assume that the condition (ii) holds. It follows from Remark
\ref{prim} that there exists $\alpha=\alpha_i \ne 0$ such that
$\alpha^M \ne 1$ for some $i$ with $1 \le i \le r$. This implies
that if we put $u=u_i$ then
$$\tr(\rho_S(u))=(\tr(\rho_0(u))^M=\alpha^M \not\in\F_2.$$
\end{proof}

Now, let us put $q'=q^2=p^{2m}$.
 We write $x \mapsto \bar{x}$ for the involution $a \mapsto a^{q}$ of $\F_q$.
 Let us consider the special unitary group $\SU_3(\F_q)$ consisting of all matrices
  $A \in \SL_3(\F_{q^2})$ which preserve a nondegenerate  Hermitian sesquilinear form on $\F_{q^2}^3$ say,
$$x,y \mapsto x_1\bar{y_3}+x_2\bar{y_2}+x_3\bar{y_1} \quad \forall
x=(x_1,x_2,x_3), y=(y_1,y_2,y_3).$$
It is well-known that the
conjugacy class of the special unitary group in $\GL_3(\F_{q^2})$
does not depend on the choice of  Hermitian form and
$\#(\SU_3(\F_q))=(q^3+1)q^3(q^2-1)$.
Clearly, for each $\beta \in \F_{q}^*$  the group $\SU_3(\F_q)$
contains the diagonal matrix $u=\diag(\beta,1,{\beta}^{-1})$ with
eigenvalues $\beta, 1, {\beta}^{-1}$; clearly, the trace of $u$ is
$\beta+{\beta}^{-1}+1$.

\begin{thm}
\label{SU3}
Suppose $G=\SU_3(\F_q)$. Suppose $V$ is an absolutely
simple nontrivial $\F_2[G]$-module. Assume that $m>1$.
If $\dim_{\F_2}(V)$ is a power of $2$ then it is equal to $q^3$.
In particular, $V$ is the Steinberg representation of $\SU_3(\F_q)$.
\end{thm}

\begin{proof}
Recall ( \cite{Gor}, p. 77, 2.8.10c), that
the adjoint representation of $G$ in $\End_{\F_{q^2}}(\F_{q^2}^3)$
splits into a direct sum of the trivial one-dimensional representation (scalars)
 and an absolutely simple $\F_{q^2}[G]$-module $\St_2$ of dimension 8 (traceless operators).
 The kernel of the natural homomorphism
$$G=\SU_3(\F_q) \to \Aut_{\F_{q^2}}(\St_2) \cong \GL_8(\F_{q^2})$$
coincides with the center $Z(G)$ which is either trivial or a
cyclic group of order $3$ depending on whether $(3,q+1)=1$ or $3$.
In both cases we get an embedding
$$G':=G/Z(G)=\U_3(q)=\PSU_3(\F_q) \subset \GL_8(\F_{q^2}).$$

If $m=2$ (i.e., $q=4$) then $G=\SU_3(\F_4)=\U_3(4)$ and one may
use Brauer character tables  \cite{AtlasB} in order to study
absolutely irreducible representations of $G$ in characteristic
$2$. Notice (\cite{AtlasB}, p. 284) that the reduction  modulo $2$
of the irrational constant b5  does not lie in $\F_2$. Using the
Table on p. 70 of \cite{AtlasB}, we conclude that there is only
one (up to an isomorphism) absolutely irreducible representation
of $G$ defined over $\F_2$ and its dimension is $64=q^3$. This
proves the assertion of the theorem in the case of $m=2,q=4$.  So
further we assume that
        $$m\ge 3,\quad q=2^m \ge 8.$$

 Clearly, for each $u \in G \subset \GL_3(\F_{q^2})$ with trace
$\delta \in \F_{q^2}$
 the image $u'$ of $u$ in $G'$ has trace $\bar{\delta}\delta-1 \in \F_q$.
 In particular,
if $u=\diag(\beta,1,\beta^{-1})$ with $\beta \in \F_q^*$  then the
trace of $u'$ is
$$t_{\beta}:=\tr(u')=(1+\beta+\beta^{-1})
(1 +\beta+\beta^{-1})-1=
 (\beta+\beta^{-1})^2.$$
 Now let us start to vary $\beta$ in the $q-2$-element set
$$\F_q\setminus\F_2=\F_q^*\setminus \{1\}.$$
One may easily check that the set of all $t_{\beta}$'s consists of
$\frac{q-2}{2}$ elements of $\F_q^*$. Since $q \ge 8$,
$$r:=\frac{q-2}{2}>\frac{q-1}{3}.$$
This implies  that $G'\subset \GL_8(\F_{q^2})$ satisfies the
conditions of Lemma \ref{trace} with $d=8$. In particular, none of
representations $\rho_S$ of $G'$ could be realized over $\F_2$
 if $S$ is a {\sl proper} subset of $\{0,1, \ldots, m-1\}$.
On the other hand, it is known (\cite{Gor}, p. 77, Example 2.8.10c) that
each absolutely irreducible representation of $G$ over $\f$ either has dimension divisible by $3$ or
is isomorphic to the representation obtained from some $\rho_S$ via
$G \to G'$. The rest is clear.
\end{proof}

\begin{thm}
\label{U3QB}
Suppose $m>1$ is an integer and let us put $q=2^m $. Let $B$ be a $(q^3+1)$-element set.
 Let $G'$ be a group acting faithfully on $B$. Assume that $G'$ contains a subgroup $G$ isomorphic to $\U_3(q)$.
 Then the $G'$-module $Q_B$ is very simple.
\end{thm}

\begin{proof}
First, $\U_3(q)$ is a simple non-abelian group,  whose order is
$q^3(q^3+1)(q^2-1)/(3,q+1)$  (\cite{Atlas}, p. XVI, Table 6;
\cite{Gor}, pp. 39--40). Second, notice that $\U_3(q)\subset G'$
acts  transitively on $B$. Indeed, the classification of subgroups
of $\U_3(q)$ (\cite{Gor}, Th. 6.5.3 and its proof, p. 329--332)
implies that each subgroup of $\U_3(q)$ has index $\ge
q^3+1=\#(B)$. This implies that $\U_3(q)$ acts transitively on
$B$. Third, we claim that this action is, in fact, doubly
transitive. Indeed,
 the stabilizer $\U_3(q)_b$ of a point $b \in B$ has index $q^3+1$ in $\U_3(q)$.
 It follows easily from the same classification that $\U_3(q)_b$  is (the image of) the
 stabilizer (in $\SU_3(\F_q)$) of a proper subspace $L$ in $\F_{q^2}^3$.
 If $L$ is a plane then counting arguments imply that the restriction of the Hermitian form to $L$
 could not be non-degenerate and therefore $\U_3(q)_b$  coincides
with  (the image of) the stabilizer of certain isotropic line
$L'\subset L\subset \F_{q^2}^3$. (The line $L'$ is the orthogonal
complement of $L$.) If $L$ is a line then counting arguments
imply that $L$ is isotropic. Hence we may always assume that
$\U_3(q)_b$  is (the image of) the stabilizer of an  isotropic
line in $\F_{q^2}^3$. Taking into account that the set of
isotropic lines in $\F_{q^2}^3$ has cardinality $q^3+1=\#(B)$, we
conclude that
 $B=\U_3(q)/\U_3(q)_b$  is isomorphic (as $\U_3(q)$-set) to the set  of isotropic lines on which $\U_3(q)$
  acts doubly transitively and we are done.

By Remark \ref{St2}, the double transitivity  implies that the $\F_2[\U_3(q)]$-module $Q_B$ is absolutely simple.
Since $\SU_3(\F_q) \to \U_3(q)$ is surjective,
the $\F_2[\SU_3(\F_q)]$-module $Q_B$ is also absolutely simple. Also, in order to prove that
$\F_2[\U_3(q)]$-module $Q_B$ is very simple, it suffices to check that the
$\F_2[\SU_3(q)]$-module $Q_B$ is very simple.

 Recall that
$\dim_{\F_2}(Q_B)=\#(B)-1=q^3=2^{3m}$.
By Theorem \ref{SU3}, there no absolutely simple
nontrivial $\F_2[\SU_3(\F_q)]$-modules, whose dimension {\sl strictly} divides $2^{3m}$.
 This implies that $Q_B$ is {\sl not} isomorphic to a tensor product
 of absolutely simple $\F_2[\SU_3(\F_q)]$-modules of dimension $>1$.
Therefore $Q_B$ is {\sl not} isomorphic to a tensor product of absolutely simple
$\F_2[\U_3(q)]$-modules of dimension $>1$.
Recall that all subgroups in $G=\U_3(q)$ different from $\U_3(q)$ itself
 have index $\ge q^3+1> q^3=\dim_{\F_2}(Q_B)$.
 It follows from Corollary \ref{Very3} that the $G$-module $Q_B$ is very simple.

\end{proof}

\section{Proof of Theorems \ref{main2}}
\label{final}
Recall that $\Gal(f) \subset \Perm(\R)$.
 It is also known that the natural homomorphism $\Gal(K) \to \Aut_{\F_2}(J(C)_2)$ factors
through the canonical surjection $\Gal(K) \twoheadrightarrow
\Gal(K(\R)/K)=\Gal(f)$ and the $\Gal(f)$-modules $J(C)_2$ and
$Q_{\R}$ are isomorphic (see, for instance, Th. 5.1 of
\cite{Zarhin2}). In particular, if the $\Gal(f)$-module $Q_{\R}$
is very simple then the $\Gal(f)$-modules $J(C)_2$ is also very
simple and therefore is absolutely simple.

\begin{lem}
\label{cor51}
If the $\Gal(f)$-module $Q_{\R}$ is very simple then
  either $\End(J(C_f))=\Z$ or $\fchar(K)>0$ and $J(C_f)$ is a supersingular abelian variety.
\end{lem}

\begin{proof}
This is Corollary 5.3 of \cite{Zarhin2}.
\end{proof}

It follows from Theorem \ref{U3QB} that under the assumptions of Theorem \ref{main2},
 the $\Gal(f)$-module $Q_{\R}$ is very simple.
 Applying Lemma \ref{cor51}, we conclude that  either $\End(J(C_f))=\Z$ or $\fchar(K)>0$ and $J(C_f)$ is
  a supersingular abelian variety.


\begin{thebibliography}{99}


\bibitem{Atlas} J. H. Conway, R. T. Curtis, S. P. Norton, R. A. Parker, R. A. Wilson, Atlas of finite groups. Clarendon Press, Oxford, 1985.

\bibitem{Curtis} Ch. W. Curtis, {\em The Steinberg character of a finite group with a $(B,\,N)$-pair}. J. Algebra {\bf 4} (1966), 433--441.

\bibitem{Gor} D. Gorenstein, R. Lyons, R. Solomon,
 The Classification of the finite simple groups, Number 3. AMS, Providence, RI, 1994.

\bibitem{Hump} J. E. Humphreys, {\em The Steinberg representation}. Bull. AMS (N.S.) {\bf 16} (1987), 247--263.


\bibitem{AtlasB} Ch. Jansen, K. Lux, R. Parker, R. Wilson,  An Atlas
of Brauer characters. Clarendon Press, Oxford, 1995.


\bibitem{Katz1} N. Katz, {\em Monodromy of families of curves:
    applications of some results of Davenport-Lewis}. In:
    S\'eminaire de Th\'eorie des Nombres, Paris 1979-80
    (ed. M.-J. Bertin); Progress in Math. {\bf 12},
    pp. 171--195, Birkh\"auser, Boston-Basel-Stuttgart,
    1981.

\bibitem{Katz2} N. Katz,  {\em Affine cohomological transforms,
    perversity, and monodromy}.
     J. Amer. Math. Soc. {\bf 6} (1993), 149--222.

\bibitem{Masser} D. Masser, {\em Specialization of
    some hyperelliptic jacobians}. In:
    Number Theory in Progress
    (eds.  K. Gy\"ory, H. Iwaniec, J. Urbanowicz), vol. I, pp. 293--307;
     de Gruyter, Berlin-New York, 1999.

\bibitem{Mori1} Sh. Mori, {\em The endomorphism rings of some abelian varieties}.  Japanese J. Math,  {\bf 2}(1976), 109--130.

\bibitem{Mori2} Sh. Mori, {\em The endomorphism rings of some abelian varieties}. II, Japanese J. Math,  {\bf 3}(1977), 105--109.

\bibitem{SerreRep} J.-P. Serre, Linear representations of finite groups, Springer-Verlag, 1977.


\bibitem{Zarhin} Yu. G. Zarhin, {\em Hyperelliptic jacobians without
complex multiplication}. Math. Res. Letters {\bf 7}(2000), 123--132.

\bibitem{Zarhin2} Yu. G. Zarhin, {\em Hyperelliptic jacobians and modular
representations}. In: Moduli of abelian varieties (C. Faber, G.
van der Geer, F. Oort, eds.), pp. 473--490, Progress in Math.,
Vol. 195, Birkh\"auser, Basel--Boston--Berlin, 2001.

\bibitem{ZarhinMRL2}  Yu. G. Zarhin, {\em Hyperelliptic jacobians without
complex multiplication in positive characteristic}. Math. Res.
Letters {\bf 8} (2001), to appear.


\end{thebibliography}
\end{document}